\def \version {2015--11--30}

\documentclass[12pt]{article}

\usepackage{amssymb}

\usepackage{tikz,ifthen}
\usetikzlibrary{calc}

\def \bsk {\bigskip}

\def \nin {\noindent}
\def \pf {\nin{\bf Proof.} \ }
\def \qed {\hfill $\Box$}

\def \cE {{\cal E}}

\def \cH {{\cal H}}
\def \vp {\varphi}

\def \NP {{\sf NP}}

\def \tmz {\begin{itemize}}
\def \etmz {\end{itemize}}

  \newtheorem{thm}{Theorem}
  \newtheorem{cor}[thm]{Corollary}
  \newtheorem{prp}[thm]{Proposition}
  
  \newtheorem{rmk}[thm]{Remark}
  \newtheorem{lem}[thm]{Lemma}
  \newtheorem{prm}[thm]{Problem}
  \newtheorem{cnj}[thm]{Conjecture}

\def \bp {\begin{prp}}
\def \ep {\end{prp}}

\def \bt {\begin{thm}}
\def \et {\end{thm}}

\def \bl {\begin{lem}}
\def \el {\end{lem}}

\def \br {\begin{rmk}}
\def \er {\end{rmk}}

\def \bcj {\begin{cnj}}
\def \ecj {\end{cnj}}

\def \bcr {\begin{cor}}
\def \ecr {\end{cor}}

\def \bpm {\begin{prm}}
\def \epm {\end{prm}}

\begin{document}

\title{\vspace{-8ex} ~~ \\
$F$-WORM colorings:\\
    Results for $2$-connected graphs
     \thanks{Research supported in part by the
      Hungarian Scientific Research Fund NKFIH/OTKA grant SNN 116095.}
  }
\author{Csilla Bujt\'as
  \quad and \quad \vspace{2ex} Zsolt Tuza\\
   \normalsize Department of Computer Science and Systems
 Technology\\
   \normalsize University of Pannonia, Veszpr\'em,  Hungary\\
   \small and\\
   \normalsize Alfr\'ed R\'enyi Institute of Mathematics\\
   \normalsize Hungarian Academy of Sciences,
     Budapest,  Hungary} \vspace{-2ex}
\date{\small Latest update on \vspace{-5ex} \version}
\maketitle

\begin{abstract}
Given two graphs $F$ and $G$, an $F$-WORM coloring of $G$ is an
assignment of colors to its vertices in such a way that no
$F$-subgraph of $G$ is monochromatic or rainbow. If $G$ has at least
one such  coloring, then it is called $F$-WORM colorable and
$W^-(G,F)$ denotes the minimum possible number of colors. Here, we
consider $F$-WORM colorings with a fixed $2$-connected graph $F$ and
prove the following three main results: $(1)$ For every natural
number $k$, there exists a graph $G$ which is $F$-WORM colorable and
$W^-(G,F)=k$; $(2)$ It is \NP-complete to decide whether a graph is
$F$-WORM colorable; $(3)$ For each $k \ge |V(F)|-1$, it is
\NP-complete to decide whether a graph $G$ satisfies $W^-(G,F) \le
k$. This remains valid on the class of $F$-WORM colorable graphs of
bounded maximum degree.  For complete graphs $F=K_n$ with $n \ge 3$
we also prove: $(4)$ For each $n \ge 3$  there exists a graph $G$
  and integers $r$  and $s$ such that  $s \ge r+2$, $G$ has $K_n$-WORM colorings with exactly $r$  and also with
 $s$ colors, but it admits no $K_n$-WORM
colorings with exactly $r+1, \dots, s-1$ colors. Moreover, the
difference $s-r$ can be arbitrarily large.



\end{abstract}

\section{Introduction}

 Given a graph $G$ and a color assignment to its vertices,
 a subgraph is   \emph{monochromatic} if its vertices have the
 same color; and
 it is \emph{rainbow} if the vertices have pairwise different colors.
 For graphs $F$ and $G$, an \emph{$F$-WORM coloring} of $G$ is an
 assignment of colors to the vertices of $G$ such that no subgraph
 isomorphic to $F$ is monochromatic or rainbow. This concept
  was introduced recently  by Goddard, Wash, and Xu \cite{GWX1}.

If $G$ has at least one $F$-WORM coloring, we say that it is
\emph{$F$-WORM colorable}. In this case,  $W^-(G,F)$ denotes the minimum
number of colors and $W^+(G,F)$ denotes the maximum number of colors
used in an $F$-WORM coloring of $G$;
 they are called the $F$-WORM \emph{lower} and \emph{upper chromatic
  number}, respectively. The
$F$-WORM \emph{feasible set} $\Phi_{_W}(G,F)$ of $G$ is the set of
those integers $s$ for which $G$ admits an $F$-WORM coloring with
exactly $s$ colors.
 Moreover, we say that  $G$   has a \emph{gap} at $k$ in
 its $F$-WORM chromatic spectrum, if $W^-(G,F)<k < W^+(G,F)$ but
  $G$ has no $F$-WORM coloring with precisely $k$ colors.
  The size of a gap is the number of consecutive integers
   missing
  from $\Phi_{_W}(G,F)$.
  If $\Phi_{_W}(G,F)$ has no gap---that is, if it contains all integers from the interval  $[W^-(G,F),W^+(G,F)]$---we say that the $F$-WORM feasible set (or the $F$-WORM chromatic
  spectrum) of $G$ is \emph{gap-free}.

The invariants
 $W^-(G,F)$ and $W^+(G,F)$ are not defined if $G$ is not $F$-WORM
 colorable.
Hence, wherever $W^-$ or $W^+$ appears in this paper,
 we assume   without further mention
  that the graph under consideration is $F$-WORM colorable.

\medskip

In the earlier works \cite{GWX,GWX1,BT-w1}, $F$-WORM colorings were considered
 for particular graphs $F$ --- cycles, complete graphs, and complete bipartite
 graphs; but mainly the cases of $F=P_3$ and $F=K_3$ were studied.
In this paper we make the first attempt towards a general theory;
 we study $F$-WORM colorings for all $2$-connected graphs $F$.
Our results presented here concern   colorability,
  lower  chromatic number, and gaps in the chromatic spectrum.

\subsection{Related coloring concepts}

A general structure class within which $F$-worm colorings can naturally
 be represented is called \emph{mixed hypergraphs}.
In our context its subclass called
 \emph{bi-hyper\-graphs} is most relevant.
 It means a pair $\cH=(X,\cE)$, where $\cE$ is a set system
 (the `edge set') over the underlying set $X$ (the `vertex set'),
 whose feasible colorings are those mappings $\vp:X\to\mathbb{N}$
 in which the set $\vp(e)$ of colors in every $e\in\cE$
 satisfies $1<|\vp(e)|<|e|$; in other words, the hyperedges are
 neither monochromatic nor rainbow.
In case of $F$-WORM colorings of a graph $G=(V,E)$ we have $X=V$,
 and a subset $e\subset V$ is a member of $\cE$ if and only if the
 subgraph induced by $e$ in $G$ contains a subgraph
 isomorphic to $F$.
For more information on mixed (and bi-) hypergraphs we recommend
 the monograph \cite{V-book}, the book chapter \cite{BTV-survey},
  and the regularly updated list of references \cite{VVweb}.

The exclusion of monochromatic or rainbow subgraphs has extensively been
 studied also separately.
Monochromatic subgraphs are the major issue of Ramsey theory,
 moreover minimal colorings fit naturally in the context of generalized
  chromatic number with respect to hereditary graph properties \cite{BBFMS},
 since the property of not containing any subgraph isomorphic to $F$
 is hereditary.

Also, forbidden polychromatic subgraphs arise in various contexts,
 most notably in a branch of Ramsey theory.
Namely, the maximum number
 of colors in an edge coloring of $G$ without a rainbow copy of $F$
 is termed anti-Ramsey number, and the number one larger --- which is
 the minimum number of colors guaranteeing a rainbow copy of $F$ in $G$
 in every coloring with that many colors --- is the rainbow number of $G$
  with respect to $F$.
We recommend \cite{FMO} for a survey of results and numerous
references.
 In particular, vertex colorings of graphs without rainbow star
$K_{1,s}$ subgraphs were studied in \cite{B+2, B+1}.

\subsection{Results}

%
%

Goddard, Wash, and Xu proved in
 \cite{GWX1} that if $G$ is $P_3$-WORM colorable,
 then $W^-(G,P_3)\le 2$.
Motivated by this, in \cite{GWX} they conjectured
 that $W^-(G,K_3)\le 2$ holds for every $K_3$-WORM colorable
 graph $G$.
 Moreover, they asked whether there is a constant $c(F)$ for
 every graph $F$ such that $W^-(G,F) \le c(F)$ for every $F$-WORM
 colorable $G$.
It is proved in~\cite{BT-w1} that the conjecture is false for
 $F=K_3$, and a finite $c(K_3)$ does not exist.
Now, we extend this result from $K_3$ to every $2$-connected graph.

 \begin{thm} \label{thm1}
 For every  2-connected graph\/ $F$ and positive integer\/ $k$,
  there exists a  graph\/ $G$ with\/ $W^-(G,F)=k$.
 \end{thm}

 What is more, the structure of those graphs is rich enough to imply
  that they are hard to recognize.
 We proved in~\cite{BT-w1} that for every $k\ge 2$ it is \NP-complete
 to decide whether the $K_3$-WORM lower chromatic number is at most $k$;
  moreover it remains hard on the graphs whose maximum degree is at most
    a suitably chosen constant $d_k$, whenever $k\ge 3$.
  It is left open whether the same is true for $k=2$.
The following general result is stronger also in the sense that for
 2-connected graphs $F$ of order $n\ge 4$ the bounded-degree
 version is available starting already from $k=n-1$ instead of $k=n$.


 \begin{thm} \label{thm2}
 For every  $2$-connected graph\/ $F$ of order $n\ge 4$ and for every integer\/ $k\ge n-1$,
  it is NP-complete to decide whether\/ $W^-(G,F) \le k$.
 This is true already on
 the class of\/ $F$-WORM colorable graphs with bounded maximum degree\/ $\Delta(G) < 2n^2$.
   \end{thm}

The decision problem of $F$-WORM colorability is proved to be \NP-complete for
 $F=P_3$ and $F=K_3$ in \cite{GWX1} and \cite{GWX}, respectively.
%
We prove the same complexity for every 2-connected $F$.

   \begin{thm} \label{thm3}
    For every   2-connected graph\/ $F$, the decision problem
  $F$-WORM colorability is NP-complete.
   \end{thm}

 Finally, we deal with the case where $F$ is a complete graph.
 We have proved in \cite{BT-w1} that there exist graphs with
 large gaps in their $K_3$-WORM chromatic spectrum. Here we show
 that this remains valid for the $K_n$-WORM spectrum with each $n \ge
 4$.
For the sake of completeness we also include the previously known case
 of $n=3$ in the formulation.

\begin{thm} \label{thm4}
    For every\/  $n \ge 3$ and\/ $\ell \ge 1$ there exist\/ $K_n$-WORM colorable graphs
   whose\/ $K_n$-WORM chromatic spectrum contains a gap of size\/~$\ell$.
          \end{thm}

In Section~\ref{sec2} we present some preliminary results and define
a basic construction. Using those lemmas, we prove
Theorems~\ref{thm1}, \ref{thm2}, and \ref{thm3} in
Section~\ref{sec3}. In Section~\ref{sec4}, we consider the case
$F\cong K_n$ and prove Theorem~\ref{thm4}.

\subsection{Standard notation}

As usual, for any graph $G$ we use the notation
 $\omega(G)$ for clique number, $\chi(G)$ for chromatic number,
  $\delta(G)$ for minimum degree, and $\Delta(G)$ for maximum degree.

\section{Preliminaries} \label{sec2}

 Here we prove a proposition on the $F$-WORM colorability and
  lower chromatic number of complete graphs; for some extremal cases
   we also consider the possible sizes of color classes.
Then, we give a basic construction which will be referred
 to in proofs of Section~\ref{sec3}.

\bp \label{lem-complete}
 For every  graph\/ $F$ of order\/ $n$, with\/ $n\ge 2$, the following hold:
 \tmz
 \item[$(i)$] For every integer\/ $s> (n-1)^2$, the complete graph\/ $K_{s}$ is not\/ $F$-WORM
 colorable.
 \item[$(ii)$] For every integer\/ $s$ satisfying\/ $1 \le s \le
 (n-1)^2$,\/ $K_s$ is\/ $F$-WORM colorable and\/ $W^-(K_s,
 F)=\left\lceil\frac{s}{n-1}\right\rceil$.
 \item[$(iii)$] In every\/ $F$-WORM coloring of the complete graph\/
 $K_{(n-1)^2}$, there are exactly\/ $n-1$ color classes each of size\/
 $n-1$.
 \item[$(iv)$] In every\/ $F$-WORM coloring of the complete graph\/
 $K_{(n-1)^2-1}$, there are exactly\/ $n-1$ color classes such that
   one of them contains\/ $n-2$ vertices
 while the other\/ $n-2$ color classes are of size\/ $n-1$ each.
 \etmz
\ep

 \pf First, observe that if $s <n$, $K_s$ contains no subgraphs isomorphic to
 $F$ and therefore, $W^-(K_s, F)=1=\left\lceil\frac{s}{n-1}\right\rceil$.
  If $s \ge n$, a subgraph isomorphic to $F$ occurs on any $n$ vertices of $K_s$.
  Hence, in an $F$-WORM coloring of $K_s$, no $n$ vertices
 have the same color and no $n$ vertices are polychromatic; on the other hand, this is also
 a sufficient condition for $F$-WORM colorability. By the
 pigeonhole principle, if  $s>(n-1)^2$, the complete graph $K_{s}$  does not have such a color
 partition, while  $K_{(n-1)^2}$ and  $K_{(n-1)^2-1}$ can be
 $F$-WORM colored only with  color classes of sizes as stated in
 $(iii)$ and $(iv)$, respectively. It also follows that for each $s\le
 (n-1)^2$, a vertex coloring of $K_s$ with $\lceil s/(n-1)\rceil$ color
 classes of size at most $n-1$ each determines an $F$-WORM
 coloring with the smallest possible number of colors.
  \qed
 \bsk

\begin{figure}[t]
\begin{center}
\begin{tikzpicture}[scale=0.6,style=thick]
\def\vr{4pt} 


 \path (-6,0) coordinate (x);
 \path (6,0) coordinate (y);
   \path (-3,2) coordinate (a1);
   \path (3,2) coordinate (a2);
    \path (-3,-2) coordinate (a3);
    \path (3,-2) coordinate (a4);
    \path (-3,-3) coordinate (a5);
    \path (3,-3) coordinate (a6);
   \path (-3.2,1.5) coordinate (b1);
   \path (-3.2,0.5) coordinate (b2);
   \path (-3.2,-0.5) coordinate (b3);
   \path (-3.2,-1.5) coordinate (b4);
   \path (3.2,1.5) coordinate (c1);
   \path (3.2,0.5) coordinate (c2);
   \path (3.2,-0.5) coordinate (c3);
   \path (3.2,-1.5) coordinate (c4);
   \path (0,0) coordinate (s1);
    \path (0,-2.5) coordinate (s2);
    \path (-6.5,0) coordinate (s3);
    \path (6.5,0) coordinate (s4);

          %
\draw (a1) -- (a2) -- (a4) -- (a3) -- (a1);
 \draw (a3) -- (a5) -- (a6) -- (a4);

 \foreach \i in {1,...,4}
{
    \draw (x) -- (b\i) [style=thin];
    \draw (y) -- (c\i) [style=thin] ;
}

    \draw (x)  [fill=white] circle (\vr);
    \draw (y)  [fill=white] circle (\vr);

  \draw(s1)node {$S$};
  \draw(s2)node {$S'$};
  \draw(s3)node {$x$};
  \draw(s4)node {$y$};
\end{tikzpicture}
\end{center}
\caption{Gadget $G_1(F)$} \label{figG1}
\end{figure}
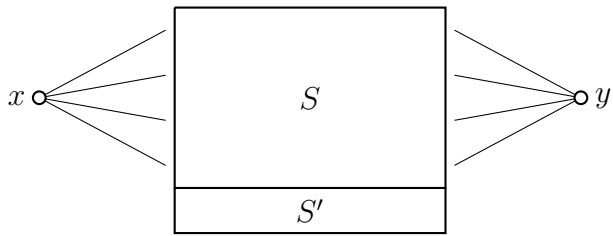

  \paragraph{Construction of the gadget  $\mathbf{G_1(F)}$.}
   For a given graph $F$ whose order is $n$ and has minimum degree $\delta\ge 2$,
   let $G_1(F)$ be the following graph. The vertex set
  is $V(G_1(F))=S\cup S'\cup\{x,y\}$ where the three sets are
   vertex-disjoint and $|S'|= n-\delta-1$,
  $|S'\cup S|=(n-1)^2-1$. Moreover, $S'\cup S$ induces a complete
  graph and the vertices $x$ and $y$ are adjacent to all vertices of
  $S$, but not to each other, neither to any vertex in $S'$.
    The vertices $x$ and $y$ will be called \emph{outer
  vertices}, while the elements of $S \cup S'$ are called \emph{inner vertices}.
  For illustration see Fig.~\ref{figG1}.

  \bl \label{lem-gadget}
   For every graph\/ $F$ of order\/ $n$ and  with minimum degree\/ $\delta
   \ge 2$, the  graph\/ $G_1(F)$ is\/ $F$-WORM colorable. Moreover, in any\/ $F$-WORM coloring of\/  $G_1(F)$,
   the outer vertices\/ $x$ and\/ $y$ get the same color which is
   repeated on exactly\/  $(n-2)$  inner vertices.
  \el
  \pf Assume that $\varphi$ is an $F$-WORM coloring  of $G_1(F)$.
  By Proposition~\ref{lem-complete}$(iv)$, $S'\cup S$ is partitioned into
  $n-1$ color classes and one of them  is of size
  $n-2$, while each further class contains exactly $n-1$ vertices.
  The color of the $(n-2)$-element color class will be denoted by
  $c^*$.

  First assume that $F \cong K_n$. Then, $S'=\emptyset$ and both
  $S\cup\{x\}$ and $S\cup \{y\}$ induce a complete subgraph on
  $(n-1)^2$ vertices.
   By Proposition~\ref{lem-complete}$(iii)$, $\varphi(x)=\varphi(y)=c^*$
   follows.

   \medskip

   If $F \ncong K_n$,  then $\delta \le n-2$ and we can take the
   following observations on $\varphi$.
  \tmz
   \item Since $S$ contains at least $n-2-|S'|=\delta-1 \ge 1$ vertices from each color
   class, we can choose an $(n-1)$-element polychromatic subset $S''$ of
   $S$. Then, on the vertex set $S'' \cup \{x\}$, which induces a complete graph, we consider a
   subgraph isomorphic to $F$. This subgraph cannot be
   polychromatic, hence the color $\varphi(x)$ (and similarly,
   $\varphi(y)$) must be assigned to at least one
   vertex of $S$.
   \item Now assume that $\varphi(x)\neq c^*$. Then, we have $n-1$
   vertices in $S' \cup S$ colored with $\varphi(x)$, and at least
   $(n-1)-|S'|= \delta$ of them are adjacent to $x$. Hence, we can
   identify a copy of $F$ monochromatic in $c^*$, in which $x$ is a
   vertex of  degree $\delta$. This cannot be the case in an
   $F$-WORM coloring. Thus, $\varphi(x)= c^*$ and similarly $\varphi(y)=
   c^*$ that proves the second part of the lemma.
   \item Consider the following coloring $\phi$ of $G_1(F)$. The color
   $c^*$ is assigned to $x$, $y$, to all vertices in $S'$, and to
   exactly $\delta-1$ vertices from $S$. The remaining $(n-2)(n-1)$
   vertices in $S$ are partitioned equally among $n-2$ further
   colors. As we used only $n-1$ colors, no subgraph isomorphic to
   $F$ can be polychromatic. Further, each color different from
   $c^*$ is assigned to only $n-1$ vertices, so no copy of $F$ can be
   monochromatic in those colors. The only color occurring on $n$
   vertices is $c^*$. But $x$ (and also $y$) shares this color with only
   $\delta -1$ of its neighbors. Therefore, we cannot have a
   subgraph  isomorphic to $F$ and monochromatic in $c^*$.
   These facts prove that $\phi$  is an $F$-WORM coloring. \qed
   \etmz
   \bsk

   \paragraph{Construction of $\mathbf{C^1(G,F, N_0)}$}
   Given an integer $N_0$, a 2-connected graph $F$ of order $n$, and a graph $G$,
   construct the following graph $C^1(G,F, N_0)$.
   If $V(G)=\{v_1, v_2,\dots , v_\ell\}$, take $N_0+1$ copies for each
   vertex $v_i$; these vertices are denoted by $v_i^0,v_i^1,\dots,v_i^{N_0}$.
   For each $1 \le i \le \ell$ and $0 \le j \le N_0-1$ take a copy
   of the gadget $G_1(F)$ such that its two outer vertices are
   identified with $v_i^j$ and $v_i^{j+1}$, respectively. The edges
   contained in these copies of  $G_1(F)$ are referred to as
   \emph{gadget-edges}.

   When we define the further edges of the construction, only the
   copy vertices of the form $v_i^{k\lceil n/2 \rceil}$ ($k \in \mathbb{N}_0$) will be used, each of them
   at most once. The sequence
   $$v_i^{0}, v_i^{\left\lceil\frac{n}{2}\right\rceil}, v_i^{2\left\lceil\frac{n}{2}\right\rceil},
   \dots, v_i^{\left\lfloor \frac{N_0}{\left\lceil\frac{n}{2}\right\rceil}\right\rfloor
   \left\lceil\frac{n}{2}\right\rceil}$$
   is called \emph{$V_i$-sequence}.

   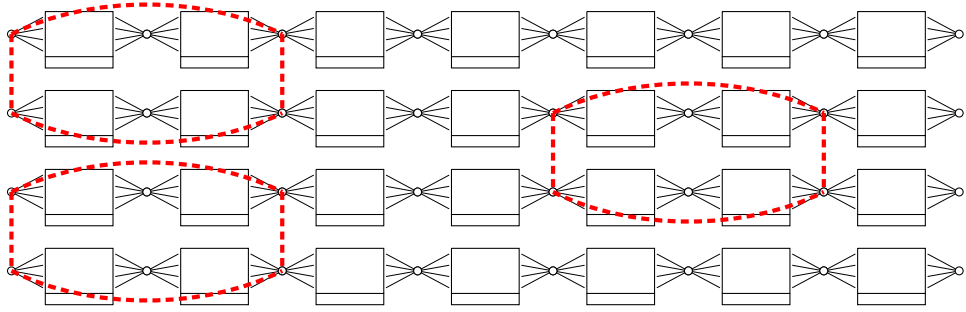
\begin{figure}[t]
\begin{center}
\begin{tikzpicture}[scale=0.15,style=thin]
\def\vr{10pt} 

   \foreach \j in {1,...,7}
   {
    \path (-6+12*\j,0) coordinate (x);
 \path (6+12*\j,0) coordinate (y);
   \path (-3+12*\j,2) coordinate (a1);
   \path (3+12*\j,2) coordinate (a2);
    \path (-3+12*\j,-2) coordinate (a3);
    \path (3+12*\j,-2) coordinate (a4);
    \path (-3+12*\j,-3) coordinate (a5);
    \path (3+12*\j,-3) coordinate (a6);
   \path (-3.2+12*\j,1.5) coordinate (b1);
   \path (-3.2+12*\j,0.5) coordinate (b2);
   \path (-3.2+12*\j,-0.5) coordinate (b3);
   \path (-3.2+12*\j,-1.5) coordinate (b4);
   \path (3.2+12*\j,1.5) coordinate (c1);
   \path (3.2+12*\j,0.5) coordinate (c2);
   \path (3.2+12*\j,-0.5) coordinate (c3);
   \path (3.2+12*\j,-1.5) coordinate (c4);
    \path (3,0) coordinate (v1);
    \path (3,-7) coordinate (v2);
    \path (3,-14) coordinate (v3);
    \path (3,-21) coordinate (v4);

          %
\draw (a1) -- (a2) -- (a4) -- (a3) -- (a1);
 \draw (a3) -- (a5) -- (a6) -- (a4);

 \foreach \i in {1,...,4}
{
    \draw (x) -- (b\i) [style= ultra thin];
    \draw (y) -- (c\i) [style=ultra thin] ;
}

    \draw (x)  [fill=white] circle (\vr);
    \draw (y)  [fill=white] circle (\vr);

  }
   \foreach \j in {1,...,7}
   {
    \path (-6+12*\j,-7) coordinate (x1);
 \path (6+12*\j,0-7) coordinate (y1);
   \path (-3+12*\j,2-7) coordinate (a11);
   \path (3+12*\j,2-7) coordinate (a21);
    \path (-3+12*\j,-2-7) coordinate (a31);
    \path (3+12*\j,-2-7) coordinate (a41);
    \path (-3+12*\j,-3-7) coordinate (a51);
    \path (3+12*\j,-3-7) coordinate (a61);
   \path (-3.2+12*\j,1.5-7) coordinate (b11);
   \path (-3.2+12*\j,0.5-7) coordinate (b21);
   \path (-3.2+12*\j,-0.5-7) coordinate (b31);
   \path (-3.2+12*\j,-1.5-7) coordinate (b41);
   \path (3.2+12*\j,1.5-7) coordinate (c11);
   \path (3.2+12*\j,0.5-7) coordinate (c21);
   \path (3.2+12*\j,-0.5-7) coordinate (c31);
   \path (3.2+12*\j,-1.5-7) coordinate (c41);

          %
\draw (a11) -- (a21) -- (a41) -- (a31) -- (a11);
 \draw (a31) -- (a51) -- (a61) -- (a41);

 \foreach \i in {11,21,31,41}
{
    \draw (x1) -- (b\i) [style= ultra thin];
    \draw (y1) -- (c\i) [style=ultra thin] ;
}

    \draw (x1)  [fill=white] circle (\vr);
    \draw (y1)  [fill=white] circle (\vr);

  }
  \foreach \j in {1,...,7}
   {
    \path (-6+12*\j,-14) coordinate (x2);
 \path (6+12*\j,-14) coordinate (y2);
   \path (-3+12*\j,2-14) coordinate (a12);
   \path (3+12*\j,2-14) coordinate (a22);
    \path (-3+12*\j,-2-14) coordinate (a32);
    \path (3+12*\j,-2-14) coordinate (a42);
    \path (-3+12*\j,-3-14) coordinate (a52);
    \path (3+12*\j,-3-14) coordinate (a62);
   \path (-3.2+12*\j,1.5-14) coordinate (b12);
   \path (-3.2+12*\j,0.5-14) coordinate (b22);
   \path (-3.2+12*\j,-0.5-14) coordinate (b32);
   \path (-3.2+12*\j,-1.5-14) coordinate (b42);
   \path (3.2+12*\j,1.5-14) coordinate (c12);
   \path (3.2+12*\j,0.5-14) coordinate (c22);
   \path (3.2+12*\j,-0.5-14) coordinate (c32);
   \path (3.2+12*\j,-1.5-14) coordinate (c42);

          %
\draw (a12) -- (a22) -- (a42) -- (a32) -- (a12);
 \draw (a32) -- (a52) -- (a62) -- (a42);

 \foreach \i in {12,22,32,42}
{
    \draw (x2) -- (b\i) [style= ultra thin];
    \draw (y2) -- (c\i) [style=ultra thin] ;
}

    \draw (x2)  [fill=white] circle (\vr);
    \draw (y2)  [fill=white] circle (\vr);

  }

  \foreach \j in {1,...,7}
   {
    \path (-6+12*\j,-21) coordinate (x3);
 \path (6+12*\j,-21) coordinate (y3);
   \path (-3+12*\j,2-21) coordinate (a13);
   \path (3+12*\j,2-21) coordinate (a23);
    \path (-3+12*\j,-2-21) coordinate (a33);
    \path (3+12*\j,-2-21) coordinate (a43);
    \path (-3+12*\j,-3-21) coordinate (a53);
    \path (3+12*\j,-3-21) coordinate (a63);
   \path (-3.2+12*\j,1.5-21) coordinate (b13);
   \path (-3.2+12*\j,0.5-21) coordinate (b23);
   \path (-3.2+12*\j,-0.5-21) coordinate (b33);
   \path (-3.2+12*\j,-1.5-21) coordinate (b43);
   \path (3.2+12*\j,1.5-21) coordinate (c13);
   \path (3.2+12*\j,0.5-21) coordinate (c23);
   \path (3.2+12*\j,-0.5-21) coordinate (c33);
   \path (3.2+12*\j,-1.5-21) coordinate (c43);

          %
\draw (a13) -- (a23) -- (a43) -- (a33) -- (a13);
 \draw (a33) -- (a53) -- (a63) -- (a43);

 \foreach \i in {13,23,33,43}
{
    \draw (x3) -- (b\i) [style= ultra thin];
    \draw (y3) -- (c\i) [style=ultra thin] ;
}

    \draw (x3)  [fill=white] circle (\vr);
    \draw (y3)  [fill=white] circle (\vr);
    }
  \draw[ultra thick, densely dashed,  color=red] (6,0) .. controls (13,3.5) and (23,3.5) ..
  (30,0);
  \draw[ultra thick, densely dashed, color=red] (6,-7) .. controls (13,-10.5) and (23,-10.5) ..
  (30,-7);
  \draw[ultra thick, densely dashed, color=red] (6,0) -- (6,-7);
  \draw[ultra thick, densely dashed, color=red] (30,0) -- (30, -7);
 \draw[ultra thick, densely dashed, color=red] (54,-7) .. controls (61,-3.5) and (71,-3.5) ..
  (78,-7);
  \draw[ultra thick, densely dashed, color=red] (54,-14) .. controls (61,-17.5) and (71,-17.5) ..
  (78,-14);
  \draw[ultra thick, densely dashed, color=red] (54,-7) -- (54,-14);
  \draw[ultra thick, densely dashed, color=red] (78,-7) -- (78,-14);
 \draw[ultra thick, densely dashed, color=red] (6,0-14) .. controls (13,3.5-14) and (23,3.5-14) ..
  (30,0-14);
  \draw[ultra thick, densely dashed, color=red] (6,-7-14) .. controls (13,-10.5-14) and (23,-10.5-14) ..
  (30,-7-14);
  \draw[ultra thick, densely dashed, color=red] (6,0-14) -- (6,-7-14);
  \draw[ultra thick, densely dashed, color=red] (30,0-14) -- (30, -7-14);

\end{tikzpicture}
\end{center}
\caption{The graph $C^1(P_4, C_4, 7)$. Supplementary edges are drawn
with dashed lines. } \label{figC1}
\end{figure}


   To finalize the construction of $C^1(G,F, N_0)$, assume $N_0 \ge (n+1)^2 \Delta(G)/4$ and consider
   the edges of $G$ one by one in an arbitrarily fixed order. When
   an edge $v_iv_j$ (with $i<j$) is treated, take the next $\lceil n/2 \rceil$
   vertices from the $V_i$-sequence and the next $\lfloor n/2
   \rfloor$ ones
   from the $V_j$-sequence, and connect them with edges to obtain an
   induced subgraph isomorphic to $F$. These edges are called
   \emph{supplementary edges}.
   For an illustration with $F=C_4$ see Fig.~\ref{figC1}.

   \bl \label{lem-C1}
   Assume that\/ $F$ is a  $2$-connected graph   of order $n$,\/ $G$ is a
   graph,  and\/ $N_0\ge \frac{(n+1)^2 \,\Delta(G)}{4}$. Then
    the graph\/ $C^1(G,F, N_0)$ satisfies the
   following  properties.
   \tmz
   \item[$(i)$] In each\/ $F$-WORM coloring\/ $\vp$ of\/ $C^1(G,F, N_0)$, the
   vertices\/ $v_i^0,v_i^1,\dots,v_i^{N_0}$ are monochromatic for each\/
   $i$ with\/    $1 \le i \le |V(G)|$. Moreover, if\/ $v_iv_j$ is an edge in\/
   $G$\/ then $\vp(v_i^0)\neq \vp(v_j^0)$.

   \item[$(ii)$] For every integer\/ $k$ with\/ $n-1 \le  k\le |V(G)|$ the graph\/    $C^1(G,F, N_0)$  is\/ $F$-WORM
   colorable with exactly\/ $k$ colors if and only if\/ $G$ is\/ $k$-colorable.
   \item[$(iii)$] For every integer\/ $k \le |V(G)|$,  there exists an\/ $F$-WORM coloring $\vp$ of\/ $C^1(G,F, N_0)$  which
   uses exactly $k$ different colors on the set of outer vertices of
   gadgets, if and only if $G$ is $k$-colorable.
   \etmz
  \el

  \pf To simplify notation, let us write $G^*=C^1(G,F, N_0)$.
   First, consider an $F$-WORM coloring $\vp$ of $G^*$.
    By Lemma~\ref{lem-gadget}, in each gadget
  $G_1(F)$  the two outer vertices have the same color. Thus, for each $i$,  the vertices
  $v_i^{0}, v_i^{1},  \dots, v_i^{N_0}$, and particularly the
  vertices contained in the $V_i$-sequence, share their color.
  We denote this color  by $\vp(V_i)$.
 By construction, if $v_iv_j$ is an edge in $G$, we have
  an $F$-subgraph in $G^*$ such that every vertex of the subgraph belongs to
  the $V_i$- or $V_j$-sequence. Since
  $F$ is not monochromatic in $\vp$, we infer that $\vp(V_i)\neq
  \vp(V_j)$. These prove $(i)$.

  Now, assume again that $\vp$ is an $F$-WORM coloring of $G^*$.
 Then, the coloring $\phi$ which assigns the color
  $\vp(V_i)$ to every vertex $v_i \in V(G)$ is a proper vertex
  coloring of $G$  and it uses precisely $|\{\vp(V_i): 1\le i\le |V(G)|\}|$ colors.
  This   proves the "only if" direction of
  $(iii)$.
 Further, we infer that    $W^-(G^*, F) \ge \chi(G)$, and if $W^-(G^*, F)\le k \le |V(G)|$ then  $G$ has a proper coloring with exactly  $k$
  colors. Since $n-1 \le W^-(G^*, F)$, this proves the "only if"
  direction of the statement $(ii)$.

  To prove the other direction, we consider an integer $k$ in the range $\chi(G)\le k \le
  |V(G)|$.
  Let $\phi$ be a  proper coloring
  of $G$ which uses the colors $1, \dots ,k$. We define a vertex coloring $\vp$ of $G^*$ as follows.
  For every $i$ and $s$, with $1 \le i \le |V(G)|$ and $0\le s \le
  N_0$, let $\vp(v_i^s):=\phi(v_i)$. Moreover, for each copy of gadget
  $G_1(F)$ whose outer vertices are $v_i^s$ and $v_i^{s+1}$, let its
  inner
  vertices be assigned with $n-1$ different colors from $1, \dots
  ,\max\{k,n-1\}$ without creating rainbow or monochromatic copies of $F$ inside the gadget.
   We can specify this assignment corresponding to Lemma~\ref{lem-gadget}. That is,
   $\vp(v_i^s)$ is repeated on all inner vertices
  nonadjacent to $v_i^s$ and on further $\delta -1$ inner
  vertices; each of the  further $n-2$ colors is assigned to exactly $n-1$
  inner vertices.

  It is clear from the definition that any $F$-subgraph  which is  contained
  entirely in one gadget or contains only supplementary edges is
  neither monochromatic nor rainbow under $\phi$. Next, we
  prove that there are no further $F$-subgraphs in $G^*$.
  First, assume that a subgraph isomorphic to $F$ contains only
  gadget edges but from at least two different gadgets. Then, this
  subgraph meets two consecutive gadgets and contains their common
  outer vertex $v_i^s$. As $s\neq 0$ and $s\neq N_0$, this outer vertex
   is a cut vertex in the subgraph
  determined by the gadget edges. Thus, $v_i^s$ would also be a cut
  vertex in the $F$-subgraph, what contradicts the
  $2$-connectivity of $F$. Therefore, such an $F$-subgraph does not
  occur in $G^*$. The only case that remains to be excluded is an
  $F$-subgraph  which contains both gadget edges and supplementary
  edges. In such a subgraph $F^*$, we would  have a vertex  which is incident to
  gadget edges and supplementary edges   as well. This vertex, say
  $v_i^r$, belongs to the $V_i$-sequence. If only the gadget edges are considered,  any further vertex
  of $V_i$ is at distance at least $n$ apart from $v_i^r$, while $F^*$ has only $n$ vertices and at least one of them
  belongs to a different $V_j$-sequence.  Hence,
  by deleting $v_i^r$ from $F^*$ we obtain a disconnected graph, one
  component of which is contained entirely in the sequence of gadgets between
  $v_i^r$ and $v_i^{r+\lceil n/2\rceil}$,
   or between $v_i^r$ and $v_i^{r-\lceil n/2\rceil}$.
   Again, this contradicts the
  $2$-connectivity of $F$. Therefore, we have only non-monochromatic and non-rainbow $F$-subgraphs, and
  $\vp$ is an $F$-WORM coloring
  of $G^*$ with exactly $k$ colors.
  This completes the proof of the
  lemma. \qed

  \section{Lower chromatic number and\\ WORM-colorability}
  \label{sec3}

  Having Lemma~\ref{lem-C1} in hand, we are now in a position to prove
  Theorems~\ref{thm1},
  \ref{thm2}, and \ref{thm3}.
 Before the proofs, we will recall the statements of the theorems.

\medskip
 \noindent
\textbf{Theorem~\ref{thm1}}. \emph{For every 2-connected graph\/ $F$
and positive integer\/ $k$,
  there exists a  graph\/ $G$ with\/ $W^-(G,F)=k$.}
  \medskip

   \pf Let $F$ be a $2$-connected graph of order $n\ge 3$.
    By Proposition~\ref{lem-complete}$(ii)$, if $1 \le k \le n-1$ and $(k-1)(n-1)<s\le k(n-1)$,
   then $W^-(K_s)=k$. Hence, we may assume $k \ge n$. We consider the graph $G^*=C^1(G,F,N_0)$
    where $G$ is a graph of chromatic number $k$,
   and $N_0= \left\lceil \frac{(n+1)^2 \,\Delta(G)}{4} \right\rceil$.
       By Lemma~\ref{lem-C1}, for every integer $k'\in [n-1, |V(G)|]$,  $G^*$ has
   an $F$-WORM coloring using exactly $k'$ colors if and only if $k' \ge
   \chi(G)$.
   Since $\chi(G)=k$ by assumption,
        this implies
   $W^-(G^*,F)=\chi(G)=k$, as desired.
   \qed

  \bsk

  \noindent \textbf{Theorem~\ref{thm2}}. \emph{For every  $2$-connected graph\/ $F$ of order $n\ge 4$ and for every integer\/ $k\ge n-1$,
  it is NP-complete to decide whether\/ $W^-(G,F) \le k$.
 This is true already on
 the class of\/ $F$-WORM colorable graphs with bounded maximum degree\/ $\Delta(G) < 2n^2$.}
  \medskip

  \pf
Let a $2$-connected graph $F$ of order $n\ge 4$ and an integer $k\ge n-1$
 be given.
 Clearly, the decision problem `\,Is $W^-(G,F) \le k$?\,' belongs
  to \NP.
To prove that it is \NP-hard (also under the assumption of
 bounded maximum degree), we apply reduction from the classical
 problem of graph $k$-colorability, which is \NP-complete
 for every $k\ge 3$.

  For a generic instance $G$ of the graph $k$-colorability
  problem, construct $G^*=C^1(G,F,N_0)$ with
   $N_0= \left\lceil \frac{(n+1)^2 \,\Delta(G)}{4} \right\rceil$.
  By Lemma~\ref{lem-C1}, $W^-(G^*,F) \le k$ if and only if $\chi(G) \le
  k$.
  Concerning the order and maximum degree of  $G^*$, we observe that
  $$|V(G^*)|=\left((n-1)^2N_0+1\right)|V(G)|$$
 and
  $$\Delta(G^*) \le \max\{(n-1)^2,
  2\left((n-1)^2-1-(n-\delta-1)\right)+\Delta(F )\} < 2n^2.
  $$
  Therefore,   the order of $G^*$ is polynomially bounded in terms
  of  $|V(G)|$ and its maximum degree satisfies the condition  given in the
  theorem. This completes the proof.
   \qed

  \bsk
   \noindent \textbf{Theorem~\ref{thm3}}. \emph{For every
      2-connected graph\/ $F$, the decision problem\/
  $F$-WORM colorability is NP-complete.}
  \medskip

 \pf Let us consider a  $2$-connected graph $F$ and  denote its order by $n$.
 The  problem is clearly in \NP. It  is proved in \cite{GWX} that the decision problem of
 $K_3$-WORM colorability is \NP-complete. Hence, we may assume that
 $n \ge 4$.
 The algorithmic
 hardness will be reduced from the decision problem of $\chi(G) \le
 n-1$ that is  \NP-complete  for each $n \ge 4$.

 For a general instance $G$ of the decision problem
 `\,$\chi(G) \le n-1$\,'
  we again begin with constructing a graph
   $C^1(G,F, N_0)$, but now with a much larger $N_0$, namely
   $$N_0= \left\lceil \frac{(n+1)^2 \,\Delta(G)}{4} \right\rceil
    + {|V(G)|-1 \choose n-1} \left\lceil \frac{n}{2} \right\rceil.
   $$
 It will be extended with further supplementary edges, as follows.

We consider those $n$-element subsets $\{i_1, \dots, i_n\}$
   of the index set $\{1, \dots, |V(G)|\}$ for which the subgraph
   induced by $\{v_{i_1}, \dots, v_{i_n}\}$ contains at least one edge.
    For each such $\{i_1, \dots, i_n\}$ we choose one vertex
   (the first one which has not been used so far)
   from each $V_i$-sequence with indices $i=i_1, \dots, i_n$,
   and take $|E(F)|$ new supplementary edges in such a way that these $n$
   vertices induce a subgraph isomorphic to $F$.
    These edges will be called supplementary edges of the second
   type. As $F$ is $2$-connected and the vertices in the $V_i$-sequences are far enough,
   this supplementation does not create any further new $F$-subgraphs different
   from the ones inserted for the selected $n$-element subsets.

Let us denote by $C^2(G,F)$ the graph obtained in this way.
  It has fewer than $|V(G)|\cdot N_0\cdot n^2$ vertices,
   which is smaller than $|V(G)|^{n+3}$ if $|V(G)|>n$.
   Therefore, once  the graph $F$ is fixed, the size of $C^2(G,F)$ is
   bounded above by a polynomial in the size of $G$.
Thus, the proof will be done if we show that $C^2(G,F)$ is
  $F$-WORM colorable if and only if $G$ has a proper vertex coloring with at
  most $n-1$ colors.

 Suppose first that $G$ admits a proper $(n-1)$-coloring $\varphi$.
This yields an $F$-WORM coloring of $C^1(G,F, N_0)$ by
 Lemma~\ref{lem-C1}, in which each $V_i$-sequence is monochromatic,
  they altogether contain precisely $n-1$ colors, and if $v_iv_j$
 is an edge in $G$ then the colors of $V_i$ and $V_j$ are different.
Then the $F$-subgraphs formed by the supplementary edges of the second
   type cannot be monochromatic, because each selected $n$-set
   $\{v_{i_1},\dots,v_{i_n}\}$ is supposed to induce at least one
   edge in $G$;
   and they cannot be rainbow $F$-subgraphs either, because only
   $n-1$ colors occur on the $V_i$-sequences.
Thus, $C^2(G,F)$ is $F$-WORM colorable in this case.

Next, assume that $\chi(G)\ge n$, and suppose for a contradiction
 that $C^2(G,F)$ admits an $F$-WORM coloring $\phi$.
Since $C^1(G,F, N_0)$ is a subgraph of $C^2(G,F)$,
 Lemma~\ref{lem-C1} implies also for the latter graph that each
  $V_i$-sequence is monochromatic in every $F$-WORM coloring,
  and any $k$-coloring of the $V_i$-sequences induced by an
  $F$-WORM coloring of $C^2(G,F)$ is a proper $k$-coloring of $G$.
Such a coloring
 necessarily uses at least $n$ colors.
Selecting an arbitrary edge $v_iv_j$ of $G$, we can extend $\{v_i,v_j\}$
 to an $n$-element set $\{v_{i_1},\dots,v_{i_n}\}$ such that
 all those vertices have mutually distinct colors.
It follows that the $F$-subgraph formed by the supplementary edges
 of the second type inserted for $\{v_{i_1},\dots,v_{i_n}\}$ is a
 rainbow copy of $F$, contradicting the assumption that $\phi$
 is an $F$-WORM coloring,


  Therefore, once $F$ is fixed according to the conditions in the theorem and $n\ge 4$,
  the decision problem of  $\chi(G) \le n-1$ can be
  polynomially reduced to the $F$-WORM colorability problem,
   and it follows that the latter problem is \NP-complete.
  \qed

  \bsk

We close this section with a positive result, implying that important
 graph classes admit efficiently solvable instances of WORM colorability.

\bp
 Let\/ $n\ge 3$ be an integer, and\/ $G$ a graph with\/ $\chi(G)=\omega(G)$.
 Then\/ $G$ is\/ $K_n$-WORM colorable if and only if\/ $\omega(G)\le(n-1)^2$.
\ep

 \pf
We know from Proposition \ref{lem-complete}$(i)$ that $K_{(n-1)^2+1}$
 is not $K_n$-WORM colorable, therefore the condition $\omega(G)\le(n-1)^2$
 is necessary.
Conversely, suppose that $\chi(G)\le(n-1)^2$.
 Take any proper coloring of $G$ with at most $(n-1)^2$ colors.
It is possible to group the color classes into exactly $n-1$ disjoint non-empty parts,
 say $C^1,\dots,C^{n-1}$, each of them consisting of at most $n-1$ colors.
(We may assume $\omega(G)\ge n$, otherwise $G$ trivially is
 $K_n$-WORM colorable.)
 Assign color $i$ to the vertices in $C^i$, for $i=1,\dots,n-1$.
Then no rainbow $K_n$ can occur because there are at most $n-1$
 colors are used, and no monochromatic $K_n$ can occur because
 each $K_n$-subgraph meets exactly $n$ color classes in the
 original proper coloring of $G$, at most $n-1$ of which
 belong to the same $C^i$.
Thus, $G$ is $K_n$-WORM colorable.
 \qed

\bsk

Since a proper coloring of a perfect graph with the minimum number
 of colors can be determined in polynomial time \cite{GLS}, we obtain:

 \bcr
  For every fixed\/ $n\ge 3$, the problem of\/ $K_n$-WORM colorability
   can be solved in polynomial time on perfect graphs.
 \ecr

  \section{Gaps in the chromatic spectrum} \label{sec4}

The following kind of graph product will play an important role in the
 proof below.
Given two graphs $G_1$ and $G_2$, the \emph{strong product}
 denoted by $G_1\boxtimes G_2$ has vertex set $V(G_1)\times V(G_2)$,
 and any two edges $u_1v_1\in E(G_1)$ and $u_2v_2\in E(G_2)$
 give rise to a copy of $K_4$ in $G_1\boxtimes G_2$ with the
 following six edges:
  $$
    \{(u_1,u_2),(u_1,v_2)\}, \quad \{(u_1,u_2),(v_1,v_2)\}, \quad \{(u_1,u_2),(v_1,u_2)\},
  $$
  $$
    \{(u_1,v_2),(v_1,u_2)\}, \quad \{(u_1,v_2),(v_1,v_2)\}, \quad \{(v_1,u_2),(v_1,v_2)\}.
  $$
Moreover, we denote by $G_1\vee G_2$ the join of $G_1$ and $G_2$,
that is
 the graph whose vertex set is the disjoint union $V(G_1)\cup V(G_2)$, and
 has the edge set
 $$
   E(G_1\vee G_2) = E(G_1)\cup E(G_2) \cup \{ v_1v_2 : v_1\in V(G_1), v_2\in V(G_2) \} .
 $$

Applying these operations, here we prove Theorem~\ref{thm4}; let us
recall its assertion.

\bsk

  \noindent \textbf{Theorem~\ref{thm4}}. \emph{For every\/  $n \ge 3$ and\/ $\ell \ge 1$
   there exist\/ $K_n$-WORM colorable graphs
   whose\/ $K_n$-WORM chromatic spectrum contains a gap of size\/~$\ell$.}

   \medskip

   \pf
 As we mentioned in the Introduction, for $K_3$
  the theorem was proved in \cite{BT-w1}.
Hence, from now on we assume $n \ge 4$.

 Consider  a triangle-free, connected graph $G$ with $\chi(G)=k\ge
 3$, and construct the  graph $G^*=(G\boxtimes K_{n-1})\vee
 K_{(n-3)(n-1)}$. When $G^*$ is obtained from $G$,
 each vertex $v_i \in V(G)$ is replaced with a complete graph on $n-1$
 vertices --- this vertex set will be denoted by $V_i$ ---
 and each edge $v_iv_j \in E(G)$ is replaced with a complete bipartite graph
 between $V_i$ and $V_j$.
  To complete the construction, we extend the graph with  $(n-3)(n-1)$ universal vertices whose set
 is denoted by $V^*$. Note that the vertex sets $V_1, \dots, V_{|V(G)|},
 V^*$ are pairwise disjoint.

 If a $K_n$ subgraph of $G^*$ meets both
  sets $V_i$ and $V_j$ (with $i\neq j$), then there exist some edges between
  these sets and hence $v_i$ and $v_j$ must be adjacent in $G$.
  Moreover, as  $G$ is triangle-free, a complete
  subgraph of $G^*$ cannot meet three different vertex sets $V_s$.
  This implies that for each $K_n$ subgraph $K$  of $G^*$ there exists an
 edge $v_iv_j\in E(G)$ such that $V(K) \subset V_i \cup V_j \cup
 V^*$. Therefore, a  vertex
 coloring $\vp$ of $G^*$ is $K_n$-WORM if and only if  the complete
 subgraph of order $(n-1)^2$  induced by $V_i \cup V_j \cup  V^*$ in
 $G^*$ is $K_n$-WORM colored for each edge $v_iv_j$ of $G$. By
 Proposition~\ref{lem-complete}, this gives the following necessary
 and sufficient condition for $\vp$ to be a $K_n$-WORM coloring:
 \tmz
 \item[$(\star)$] For each $v_iv_j \in E(G)$, $\vp$ uses exactly
 $n-1$ colors on $V_i \cup V_j \cup  V^*$, and each color occurs on
 exactly $n-1$ vertices of this complete subgraph.
 \etmz

 Now, we assume that $\vp$ is a $K_n$-WORM coloring of $G^*$. We make the following observations.
 \tmz
 \item Since
 there exist $K_{(n-1)^2}$-subgraphs, $\vp$ uses at least $n-1$
 colors. On the other hand, by $(\star)$ a $K_n$-WORM coloring is
 obtained if each of the colors $1,2, \dots, n-1$ occurs on exactly $n-3$ vertices from
 $V^*$, and on exactly one vertex from each $V_i$. This proves
 $W^-(G^*,K_n)=n-1$.
 \item If $\vp$ uses exactly $n-1$ colors on $V^*$, it follows from  $(\star)$
  that no further colors appear on the sets $V_i$.
 \item If $|\vp(V^*)|=n-2$, then for each $v_iv_j \in E(G)$ the set
 $\vp(V_i \cup V_j)$ contains exactly one color different from those
 in $\vp(V^*)$. We have two cases. If there exists a monochromatic
 $V_s$, its color $c^*$  appears on $n-1$ vertices in $V_s$.
 By  $(\star)$,  $c^* \notin \vp(V^*)$ follows, and also that for every
 neighbor $v_p$ of $v_s$, $c^* \notin \vp(V_p)$. Then, $|\vp(V^* \cup
 V_p)|=n-2$ and for each neighbor $v_q$ of $v_p$, the vertex set
 $V_q$ in $G^*$ must be monochromatic in a color not included in
 $\vp(V^*)$. As $G$ is connected, this property propagates  along
 the edges and for every adjacent vertex pair $v_i,v_j$, one of the sets
 $V_i$ and $V_j$ is monochromatic and the other is not. This
 gives a bipartition of $G$, which contradicts our assumption
 $\chi(G) \ge 3$. In the other case, there is no monochromatic
 $V_i$,   therefore the $n-1$ vertices of the $(n-1)$st color
  of $V_i \cup V_j \cup  V^*$ have to be distributed between
   $V_i$ and $V_j$.
  This implies
 $$\vp(V_i   \cup  V^*)= \vp(V_i \cup V_j \cup  V^*)= \vp( V_j \cup
 V^*)$$
 for every pair $i,j$ with $v_iv_j \in E(G)$. By the
 connectivity of $G$, we conclude that $|\vp(G^*)|=n-1$.
 \item Assume that $|\vp(V^*)|=n-3$. Then, each of these $n-3$ colors occurs on exactly
 $n-1$ vertices of $V^*$ and occurs on no further vertices of $G^*$.
 Moreover, for each $v_iv_j \in E(G)$, the vertices in
 $V_i \cup V_j$ are colored with exactly two colors such that each
 color is assigned to exactly $n-1$ vertices.
 If there is a non-monochromatic $V_s$, then $\vp(V_s)=\vp(V_p)$ for
 every $p$ satisfying $v_sv_p \in E(G)$. Then, since $G$ is connected,
 this  equality will be also valid if $v_sv_p \notin E(G)$.
 Therefore,  we have only $n-1$ different colors on the vertices of $V(G^*)$,
 again. On the other hand,
  if every $V_i$ is made monochromatic by $\vp$, the condition
 $(\star)$ is satisfied if and only if $(i)$ the color of $V_i$ is not in
 $\vp(V^*)$; and $(ii)$ for every adjacent vertex pair $v_i,v_j$ of
 $G$, the colors $\vp(V_i)$ and $\vp(V_j)$ are disjoint.
 Conditions $(i)$ and $(ii)$ imply that the color assignment $\phi$
 defined as $\phi(v_i)=\vp(V_i)$ gives a proper vertex coloring of
 $G$ with $|\vp(V(G^*))|-n+3$ colors. Hence, this type of $K_n$-WORM
 coloring of $G^*$ can be constructed such that the number of used
 colors is one from the range $\chi(G)+n-3, \dots , |V(G)|+n-3$.
\etmz

We have proved that the $K_n$-WORM feasible set of $G^*$ is
$$\{n-1\} \cup \{k+n-3, \dots , |V(G)|+n-3\} .
$$
 If we choose a triangle-free connected graph $G$ with
 $\chi(G)=k= \ell+3$, the gap in the feasible set $\Phi_W(G^*,K_n)$ is
  of size $\ell$. \qed

\end{document}